\documentclass[conference]{worldcomp}

\usepackage[hmargin=.75in,vmargin=1in]{geometry}
\usepackage[american]{babel}
\usepackage[T1]{fontenc}
\usepackage{times}
\usepackage{caption}

\usepackage{textcomp}
\usepackage{epsfig,graphicx}
\usepackage{xcolor}
\usepackage{amsfonts,amsmath,amssymb}
\usepackage{fixltx2e} % Fixing numbering problem when using figure/table* 
\usepackage{booktabs}

\title{\bf Wentzel-Freidlin estimates for jump processes in semi-group theory: upper bound}           %%%% Replace with your title.

\author{
{\bfseries R. L\'eandre}\\
Institut de Math\'ematiques. Universit\'e de Bourgogne, Dijon, 21000, France\\
}

\begin{document}

\maketitle                        %%%% To set Title and Author names.

\begin{abstract}%%%% Replace with your abstract.

We give a transtation in semi-group theory of Wentzel-Freidlin estimates for Poisson process.
We consider the case of the upper bound.
\end{abstract}

\vspace{1em}
\noindent\textbf{Keywords:}
 {\small  Large deviations. Semi-groups.} %%%% Replace with your keywords

%%%%%%%%%%%%%%%%%%%%%%%%%%%%%%%%%%%%%%%%%%%%%%%%%%%%%%%%%%%%

\section{Introduction}
The object of  the large deviation theory [1] is to estimate the logarithm of the probability of rare events. Bismut [2] pointed out the relationship
 between large deviation estimates and the Malliavin Calculus in order to get short time asymptotics of heat kernels associated to diffusion semi-groups.
This relationship was fully performed by L\'eandre in [3], [4].
 The reader interested in short time asymptotics of heat-kernels by using the Malliavin Calculus as a tool
 can look at the 
review of L\'eandre [5], Kusuoka [6], and Watanabe [7].

L\'eandre  has translated plenty of tools of stochastic analysis in semi-group theory  by using the fact that in the proof of these tools,
 there are suitable stochastic
 differential equations and therefore suitable parabolic equations which appear. We refer to the review of L\'eandre [8], [9] on that.

 We consider the symbol $(x,\xi) \rightarrow H(x,\xi)$ of the operator. We consider the operator $L^h$ associated to the symbol $H^h(x, \xi) = 
H(x,h\xi)$. We consider the heat-equation associated to $1/hL^h$ and his behaviour when $h \rightarrow 0$. Since in this case,
 we have a probability semi-group, it is convenient to study lower-bound and upper-bound
of the behaviour when $h\rightarrow 0$ of its solution. In semi-classical analysis [10], people consider the Schroedinger equation associated to 
$1/hL^h$. In such a case there is no representation of the solution of it by using probability measures, and people look asymptotics expansion of the solution of the 
Schroedinger equation whose main term is oscillatory.

\section{Statement of the main theorem}

Let $\mu(x,dz)= g(x,z)dz$ be a positive measure (it is the Levy measure associated to the jump process) on $\mathbb{R}^d$ such that $\sup_x\int_{\mathbb{R}^d}\vert z \vert^2\mu(x,dz) < \infty$.We suppose that $(x,z) \rightarrow g(x,z)$  
is continuous if $z \not= 0$.
 We introduce the Hamiltonian defined for $(x, \xi) \in \mathbb{R}\times \mathbb{R}^d$
\begin{equation}H(x, \xi) = \int_{\mathbb{R}^d}(\exp[<z,\xi>]-1-<z,\xi>)\mu_x(dz)\end{equation}

{\bf{Hypothesis H.1}}$H$ is continuous, bounded uniformy in x by $H_1(\xi)$, a convex function on $\mathbb{R}^d$.

We consider the Legendre transform of $H$
\begin{equation}L(x, \alpha) = \sup_{\xi}\{<\alpha,\xi>-H(x,\xi)\}\end{equation}

{\bf{Hypothesis H.2:}}The function $(x, \alpha) \rightarrow L(x, \alpha)$ is finite, twice differentiable in $\alpha$. For all $R$, there exists $M$ and $m$ such that 
$L(x, \alpha) \leq M$, $\vert {\partial \over \partial \alpha}L(x, \alpha)\vert \leq M$ and such that
 the Hessian in $\alpha$ of $L$ is larger than $mI_d$ if $\vert \alpha\vert \leq R$. Moreover
\begin{equation}\sup_{\vert x-x'\vert<\delta'}\sup_\alpha{\vert L(x',\alpha)-L(x,\alpha)\vert \over 1+L(x,\alpha)} \rightarrow 0\end{equation}
when $\delta' \rightarrow 0$.

{\bf{Hypothesis H.3:}} The Legendre transform $L_1$ of $H_1$ satisfies for all $C$ to, $L_1(\alpha) \geq C \vert \alpha \vert$ for $\vert \alpha \vert \geq K_C$.

We consider a piecewise $C^1$ curve $\phi(t)$ and we consider the action:
\begin{equation} S(\phi) = \int_0^1L(\phi(t), d/dt\phi(t))dt\end{equation}
We put
\begin{equation}l(x,y) = \inf_{\phi(0)=x, \phi(1) = y}S(\phi) \end{equation}
Under the previous assumption, $(x,y) \rightarrow l(x,y)$ is continuous.
We define the generator $L^h$ defined on smooth functions with bounded derivatives at each order
\begin{multline} L^hf(x) = \\ \int_{\mathbb{R}^d}(f(x+hz)-f(x)-h<z,f'(x)>)\mu_x(dz)\end{multline}

Under these previous assumptions, we will get in the sequel Markovian semi-group. In particular, $1/hL^h$ generates a semi-group $P_t^h$.
In such a case, $P_t^h$ is a semi-group in probability measures ([11], [1])

{\bf{Theorem}} {\it{(Wentzel-Freidlin [1])Let}} $O$ {\it{ be an open ball of}} $\mathbb{R}^d$. {\it{When}} $h \rightarrow 0$,
\begin{equation}\overline{\lim} h\,LogP_1^h[1_O](x) \leq -\inf_{y \in O}l(x,y)\end{equation}

In [12], we have proved the opposite inequality. For the case of diffusion, we refer to [13], [14].

{\bf{Remark:}}$H(x,i\xi)$ is the symbol associated to $L^1$ ([11]).

\section{Proof of the theorem}
It is the translation in semi-group theory of the proof of [1], p 150-151-152-153.

Let us consider $C \in \mathbb{R}^d$. We consider the generator on $\mathbb{R}^d \times \mathbb{R}$
\begin{multline}\overline{L}^hf(x,y) =\\  {1\over h} \int_{\mathbb{R}^d}(f(x+hz,y)-f(x,y)-h<z,f'_x(x,y)>)\mu_x(dz) +\\
{1\over h}f'_y(x,y)H(x,C)\end{multline}
We remark that
\begin{equation}\overline{L}^h\exp[<C,x>-y] = 0 \end{equation}
$\overline{L}^h$ generates a semi-group $\overline{P}^h_t$. By an argument similar to [15],, we deduce that:
\begin{equation} \overline{P}^h_t[\exp[<C,.>-(.)]](x,0) = \exp[<C,x>]\end{equation}
It is the analog in semi-group theory of the celebrated exponential martingales of  stochastic Calculus. But $H(x,C) \leq H_1(C)$ such that
\begin{equation} 1 \geq P_t^h[\exp[<C,.-x>]](x) \exp[-t/hH_1(C)]\end{equation}
We deduce that if $B(x,R)$ is the ball of center $x$ and radius $R$
\begin{equation}P_t^h[1_{B(x,R)^c}](x) \leq \exp[t/hH_1(C)]\exp[-\vert C \vert R/h]\end{equation}
We deduce that if $R$ is enough big, that $P_t^h[1_{B(x,R)^c}](x)$ has a quick exponential decay. We can replace the inequality (12) by
\begin{multline}P_t^h[1_{B(x,R)^c}](x)\\ \leq \sum_i \exp[t/hH_1(C)]\exp[-<C,R_i>/h]\end{multline}
where $R_i$ is the a finite set of element of $\mathbb{R}^d$ such that $\vert C \vert R = <C, R_i>$ for one of these. Therefore
\begin{equation}P_t^h[1_{B(x,R)^c}](x) \leq \sum_i \exp[-t/hL_1(R_i/t)]\end{equation}
By {\bf {Hypothesis H.3}}, we deduce that if $t$ is small $P_t^h[1_{B(x,R)^c}](x)$ has a quick exponential decay if $R$ is very small.

We consider $\Delta t$ small and take $t_k = k \Delta t$. Because we consider a semi-group in probability measures, we deduce on the set of 
polygonal curves $l$  $((0,x), (t_1,x_1), ..., (t_n, x_n))$ $(t_n = 1)$ a probability measure $W_{\Delta t}^h(x)$. The problem is to get an upper-bound of
$W_{\Delta t}^h(x)[1_O(x_n)]$. By using the previous exponential inequalities,, we can consider only the polygonal curves $l$ such that $\vert x_i \vert \leq R$ and such that if $\Delta t$ 
is chosen small enough that $\vert x_{i+1}-x_i\vert \leq \delta$ for a very small $\delta$. We get a set $E$ of polygonal curves.

\begin{multline}W_{\Delta t}^h(x)[1_E 1_O(x_n)]\\ \leq W_{\Delta t}^h(x)[1_E \exp[S(l)/h]] \exp[-\inf_{y \in O}l(x,y)/h]\end{multline}
By our assumption of the choice of $l \in E$, we have
\begin{multline} S(l) = \sum_{i=0}^{n-1}\int_{t_i}^{t_{i+1}}L(l_s, {x_{i+1}-x_i \over \Delta t})ds \leq \\
(1+ \chi)\Delta t \sum_{i=0}^{n-1}L(x_i, {x_{i+1}-x_i \over \Delta t}) + \chi =\\ \sum\tilde{S}(x_i,x_{i+1}) + \chi\end{multline}
for a small $\chi$.

The sequel follows closely the lines of [1] p 152. We can choose some points $\alpha_i$ such that if we put
\begin{equation}L'(x, \alpha) = \sup_i [L(x, \alpha_i) + <{\partial \over \partial \alpha}L(x, \alpha_i), \alpha-\alpha_i>]\end{equation} 
then
\begin{equation}L(x, \alpha)-L'(x, \alpha) \leq \chi\end{equation}
for $\vert \alpha \vert \leq R$. We can recognize in $L'(x, \alpha)$ the expression
\begin{equation}L'(x, \alpha) = \sup_i[<\beta_i, \alpha>-H(x, \beta_i)]\end{equation}
where $\beta_i = {\partial \over \partial \alpha}L(x, \alpha_i)$ depends on $x$. It is enough to show by the Markov property that
\begin{equation}\sup_xP_{\Delta t}^h[\exp[\tilde{S}(x,.)/h]](x)\end{equation}
has only a small exponential blowing up when $h \rightarrow 0$. By using (18), it is enough to show that
\begin{equation}\sup_xP_{\Delta t}^h[\Delta t (1+\chi)L'(x,(.-x)/\Delta t)/h]](x)\end{equation}
has a small exponential blowing up when $h \rightarrow 0$. By using the inequality
\begin{equation}\exp[\sup a_i] \leq \sum \exp[a_i])\end{equation} it is enough to show that
\begin{equation}\sup_xP_{\Delta t}^h[\exp[(1+\chi)/h[<\beta_i,.-x>-\Delta tH(x, \beta_i)]](x)\end{equation}
has a small exponential blowing-up when $h \rightarrow 0$.
We consider the semi group $\overline{P}_t^{h,i,y_1}$ on $\mathbb{R}^d \times \mathbb{R}$ associated to 
\begin{multline}\overline{L}^{h,i,y_1} =\\  {1\over h} \int_{\mathbb{R}^d}(f(x+hz,y)-f(x,y)-h<z,f'_x(x,y)>)\mu_x(dz) +\\
{1\over h}f'_y(x,y)H(x,(1+\chi)\beta_i(y_1))\end{multline}
because $\beta_i$ depends on $y_1$.
We have as in (8)
\begin{multline}\overline{P}^{h,i,x}_t[\exp[<(1+\chi) \beta_i,.>-(.)]](x,0) =\\ \exp[<(1+\chi)\beta_i,x>]\end{multline}
Since we have a semi-group in probability measures, it is enough to estimate by the Hoelder  inequality
\begin{equation}\overline{P}_t^{h,i,x}[\exp[y-tH(x, \beta_i)(1+\chi)/h]](x,0) = \gamma_t\end{equation}
We don't enter in the technicalities which arise from the change of exponent when we use the Hoelder inequality.
\begin{multline}d/dt \gamma_t = \\ \overline{P}_t^{h,i,x}[\exp[y-tH(x, \beta_i)(1+\chi)/h]\\1/h(H(.,(1+\chi)\beta_i)
-H(x,\beta_i)(1+\chi))](x,0) \end{multline}
We distinguish if $\vert H(.,(1+\chi)\beta_i)
-H(x,\beta_i)(1+\chi)\vert \geq \delta'$ or not and we remark by argument of the beginning of the proof for $t$ small enough, for a small $\delta'$ and a big $C$
\begin{multline} \overline{P}_t^{h,i,x}[\vert H(.,(1+\chi)\beta_i)
-H(x,\beta_i)(1+\chi)\vert \geq \delta'](x,0)\\ \leq \exp[-C/h]\end{multline}
For this same probability law, $y$ and $tH(x, \beta_i)(1+\chi)/h$ are bouded by $C_1t/h$. We deduce that for a small $\delta'$, we get 
the inequality for a big $C$
\begin{equation}\vert d/dt \gamma_t \vert \leq \delta'/h\gamma_t + \exp[-C/h]\end{equation} The result arises by Gronwall lemma.$\diamondsuit$

\section{Conclusion}We have translated in this note as well as in [12] some basical tools of the stochastic analysis on Poisson processes. Others basical tools were
translated in [15] and [16].

 %%%%%%%%%%%%%%%%%%%%%%%%%%%%%%%%%%%%%%%%%%
%%
%% Reference
%% Below is an example of bibliography that contains all entries within this document.
%% You can also let BibTeX generate your bibliography by inserting the following two commands:
%%
%% \bibliographystyle{IEEEtran}
%% \bibliography{<your_bibliography_file_1>,<your_bibliography_file_2>,...}

\begin{thebibliography}{1}
\providecommand{\url}[1]{#1}
\csname url@rmstyle\endcsname
\providecommand{\newblock}{\relax}
\providecommand{\bibinfo}[2]{#2}
\providecommand\BIBentrySTDinterwordspacing{\spaceskip=0pt\relax}
\providecommand\BIBentryALTinterwordstretchfactor{4}
\providecommand\BIBentryALTinterwordspacing{\spaceskip=\fontdimen2\font plus
\BIBentryALTinterwordstretchfactor\fontdimen3\font minus
  \fontdimen4\font\relax}
\providecommand\BIBforeignlanguage[2]{{%
\expandafter\ifx\csname l@#1\endcsname\relax
\typeout{** WARNING: IEEEtran.bst: No hyphenation pattern has been}%
\typeout{** loaded for the language `#1'. Using the pattern for}%
\typeout{** the default language instead.}%
\else
\language=\csname l@#1\endcsname
\fi
#2}}

\bibitem{1}A.D. Wentzel and M.J. Freidlin,
\newblock{\em " Random perturbations of dynamical systems"}
\newblock Springer, 1984.

\bibitem{2}J.M. Bismut,
\newblock{\em "Large deviations and the Malliavin Calculus"}
\newblock. Progress in Maths 45, Birkhauser, 1984.

\bibitem{3}R. L\'eandre,
"Estimation en temps petit de la densit\'e d'une diffusion hypoelliptique," \emph{ C.R.A.S. S\'erie I}, vol. 301, pp. 801--804, 1985.

\bibitem{4}R. L\'eandre,
" Majoration en temps petit de la densit\'e d'une diffusion d\'eg\'er\'ee,"  \emph{P.T.R.F.}, vol. 74, pp. 289--294, 1987.

\bibitem{5}R. L\'eandre,
" Applications quantitatives et qualitatives du Calcul de Malliavin," in \emph {French Japanese Seminar}, 1988, p. 109. 
L.N.M.  1322, Springer. English translation in \emph{Geometry of Random motion}, 1988, p. 173, 
 Contemporary Maths 73, A.M.S.

\bibitem{6}S. Kusuoka, "More recent theory of the Malliavin Calculus," \emph{Sugaku Expositions}, vol. 5, pp. 155--173, 1992.

\bibitem{7}S. Watanabe, "Stochastic analysis and its applications," \emph{Sugaku Expositions}, vol.  5, pp. 51--69, 1992.

\bibitem{8}R. L\'eandre, " Applications of the Malliavin Calculus of Bismut type without probability",  in \emph{ Simulation,
Modelling and Optimization 2006}, C.D. W.S.E.A.S.   2006, p. 559.
\emph{WSEAS transactions on mathematics}, vol. 5, pp. 1205--1211, 2006.

\bibitem {9}R. L\'eandre, " Malliavin Calculus of Bismut type in semi-group theory,"  \emph{Far East Journal of Mathematical Sciences}, 
vol. 30, pp. 1--26, 2008.

\bibitem{10}V.P. Maslov and M.V. Fedoriuk,
\newblock{\em " Semiclassical approximation in quantum mechanics."}
\newblock Reidel, 1981.

\bibitem{11}N. Jacob,
\newblock{\em "Pseudo differential operators and Markov processes. II. Generators and their potential theory."}
\newblock Imperial College Press,  2002.

\bibitem{12}R. L\'eandre, " Wentzel-Freidlin estimates for jump process in semi-group theory: lower bound" in \emph{Int. Conf. Dif. Geometry.
Dynamical. System}, \emph{B.S.G. Proceeding}, vol 17, pp 107--113, 2010.

\bibitem{13}R. L\'eandre, "Wentzel-Freidlin estimates in semi-group theory" in \emph{Control, Automation, Robotics and Vision},
 C.D., I.E.E.E.  2008, p. 223.

\bibitem{14}R. L\'eandre, "Varadhan estimates without probability: upper bound," \emph{W.S.E.A.S transactions of mathematics},
 vol. 7, pp. 244--253, 2008.

\bibitem{15}R. L\'eandre, "Malliavin calculus of Bismut type for Poisson processes without probability",   in \emph{Fractional order systems},
\emph{ J.E.S.A.}, vol.  42, pp. 715--733, 2008.

\bibitem{16}R. L\'eandre, "Girsanov transformation for Poisson processes in semi-group theory" in \emph{Num. Ana. Appl. Mathematics}, 2007, p. 336
A.IP. Proceedings 936, Amer. Inst. Phys.

\end{thebibliography}
%%
%% Note that you need to make sure that LaTeX (BibTeX) can find IEEEtrans.bst in your system.
%% If you are unsure about that, just place IEEEtrans.bst in the same directory where your LaTeX source files reside.
%%
%%%%%%%%%%%%%%%%%%%%%%%%%%%%%%%%%%%%%%%%%%%%%%%%%%%%%%%%%%%%%
%%% Below thebibliography environment will be automatically created in a different file (your_file_name.bbl) 
%%% if you use BibTeX and specify IEEEtrans.bst.

\end{document}